# ERGODIC THEORY FOR SDES WITH EXTRINSIC MEMORY[1]

By M. Hairer and A. Ohashi

*University of Warwick and Universidade Estadual de Campinas*

We develop a theory of ergodicity for a class of random dynamical systems where the driving noise is not white. The two main tools of our analysis are the *strong Feller property* and *topological irreducibility*, introduced in this work for a class of non-Markovian systems. They allow us to obtain a criteria for ergodicity which is similar in nature to the Doob–Khas'minskii theorem.

The second part of this article shows how it is possible to apply these results to the case of stochastic differential equations driven by fractional Brownian motion. It follows that under a nondegeneracy condition on the noise, such equations admit a unique adapted stationary solution.

**1. Introduction.** Ergodic properties of Markovian systems have been intensively studied, especially in the context of stochastic differential equations (henceforth abbreviated as SDEs). Many authors have been studying the problem of ergodicity for Markovian systems induced by finite- and infinite-dimensional stochastic equations driven by a Brownian motion. A good summary of the current state of research in this area can be found in the monographs [4, 8, 13]. The asymptotic behavior of processes driven by a noise with nontrivial time correlations seems to be much less well understood, although substantial progress has been made in the framework of the theory of random dynamical systems [2]. However, framework takes a rather "deterministic" approach and is mainly suitable for the study of the random equivalent of the objects from the theory of ordinary dynamical systems. A natural question is whether one can take a more "probabilistic" approach and obtain statements that are similar in spirit and in scope to the ones obtained in the Markovian case. This is the program that we start to carry out in this work. Our main goal is to obtain a criterion for the existence and

Received March 2006; revised March 2006.

[1]Supported by CAPES-Brasil Grant BEX-1534041.

*AMS 2000 subject classifications.* 60H10, 60G10, 26A33.

*Key words and phrases.* Non-Markovian processes, ergodicity, fractional Brownian motion.

uniqueness of an "invariant measure" (in a sense to be made precise) that are comparable in scope to the existing criteria for Markov processes.

More precisely, we are interested in providing a generalization of the widely used result attributed to Doob and Khas'minskii which states that a Markov process which is strong Feller and topologically irreducible can have at most one invariant measure (see, e.g., [4], Proposition 4.1.1 and Theorem 4.2.1). The obvious question that arises is how to formulate a useful generalization of the strong Feller property in non-Markovian situations. This question will be answered, to a certain extent, in the framework of "stochastic dynamical systems" (SDS) developed in [7]. Roughly speaking, an SDS is simply a random dynamical system which is reformulated in such a way that one sees how new randomness comes into the system as time evolves. One characteristic of this point-of-view is that it automatically discards invariant measures that are not measurable with respect to the past; see [2] for this terminology. Note that this is actually a desirable feature if one wishes to obtain a natural generalization of the concept of "invariant measure" from the theory of Markov processes. For example, in the case of a diffusion on the circle with a nontrivial drift, the theory of Markov processes yields the existence of a unique invariant measure. The theory of random dynamical systems, on the other hand, yields *two* distinct invariant measures, but one of them is measurable with respect to the future and corresponds to an unstable random fixed point. Even though such invariant measures correspond to stationary solutions of the corresponding SDE, they are "unphysical" in the sense that they can only be realized by preparing the initial condition in a state that depends on the whole future of the driving noise. The main result of this first, "abstract," part of the present article is Theorem 3.10 below.

As a test of the relevance of our criteria, we then show that it can be applied to the case of SDEs driven by fractional Brownian motion (fBm). The choice of fractional noise as driving noise (rather than, e.g., an Ornstein–Uhlenbeck process) is motivated by the following arguments:

1. one cannot reduce it to a Markovian situation without adding infinitely many degrees of freedom;
2. it presents long range correlations and therefore does not reduce to white noise in the limit of large time rescalings;
3. it is very well studied, so that many a priori estimates are available in the existing literature;
4. it appears naturally as the only continuous scale-invariant Gaussian process.

This article should be considered as a sequel to the work [7], where SDEs driven by *additive* fractional noise were considered. In this situation, a coupling argument allowed it to be shown that such SDEs possess a unique invariant measure under quite general conditions. Unfortunately, this argument presented two major drawbacks. First, it was very difficult to follow and hard to analyze because of the long-range correlations of the driving noise. Second, the coupling construction used the additivity of the noise in an essential way, making the argument unsuitable to the study of equations driven by multiplicative noise.

In this work, we consider equations driven by nondegenerate multiplicative noise, that is, we study the SDE

$$(1.1) \qquad dx_t = f(x_t)\,dt + \sigma(x_t)\,dB_H(t), \qquad x(0) = x_0 \in \mathbf{R}^d,$$

where $f : \mathbf{R}^d \to \mathbf{R}^d$, $\sigma : \mathbf{R}^d \to \mathrm{M}_{d \times d}$ (where $\mathrm{M}_{d \times d}$ denotes the space of $d \times d$ matrices) and $B_H$ is a $d$-dimensional fractional Brownian motion with Hurst parameter $H$. In other words, it is a centered $d$-dimensional Gaussian process with continuous sample paths, $B_H(0) = 0$ and covariance

$$\mathbf{E}(B_H^i(t) - B_H^i(s))(B_H^j(t) - B_H^j(s)) = \delta_{ij}|t - s|^{2H}$$

for $t, s \in \mathbf{R}$ and $i, j = 1, \ldots, d$. We will assume throughout this work that $H$ is strictly greater than $1/2$ so that the integral appearing in the right-hand side of (1.1) may be considered pathwise as a Riemann–Stieltjes integral. We believe that this restriction could be weakened by considering noise spaces of "rough path" type (see, e.g., [6, 12]), but this would raise additional difficulties that we do not wish to address here. A pair $(x, B_H)$ of continuous stochastic processes is called a *solution* to (1.1) if $B_H$ is a fBm and the integrated form of (1.1) holds almost surely for all times. We call such a solution *adapted* if for every $t$, $x(t)$ and $\{B_H(s)\}_{s \geq t}$ are conditionally independent, given $\{B_H(s)\}_{s \leq t}$.

In order to ensure the global existence of solutions and in order to have some control over it, we make, for most of this paper, the following assumptions on the coefficients $f$ and $\sigma$.

(H1) *Regularity*: Both $f$ and $\sigma$ are $\mathcal{C}^\infty$. Furthermore, the diffusion coefficient $\sigma$ and the derivatives of $f$ and $\sigma$ are globally bounded:

$$(1.2) \qquad \sup_{x \in \mathbf{R}^d} (|\sigma(x)| + |Df(x)| + |D\sigma(x)| < \infty).$$

(H2) *Nondegeneracy*: $\sigma(x) \in \mathrm{M}_{d \times d}$ is invertible and $\sup_{x \in \mathbf{R}^d} |\sigma^{-1}(x)| < \infty$.

(H3) *Dissipativity*: There exists $C > 0$ such that

$$\langle f(x), x \rangle \leq C(1 - \|x\|^2) \qquad \forall x \in \mathbf{R}^d.$$

The main result of this paper is the following.

THEOREM 1.1. *If the coefficients of the SDE* (1.1) *satisfy assumptions* (H1)–(H3), *then it has exactly one adapted stationary solution.*

The remainder of this paper is organized in the following way. After fixing the notation and recalling some results from [7] in Section 2, we formulate and state the main abstract result in Section 3. Section 4 is devoted to ensuring that the abstract framework constructed in Section 2 can be applied to the SDE (1.1). It also provides the a priori bounds required to ensure the existence of an invariant measure for such systems. We then spend most of Section 5 proving that the generalization of the strong Feller property formulated in Section 3 does indeed hold for (1.1). This allows us to obtain Theorem 1.1 as a simple corollary.

**2. Preliminaries.** In this section, we fix the basic notation that we use in this paper and recall some basic definitions and results from [7]. Given a product space $\mathcal{X} \times \mathcal{Y}$, we denote by $\Pi_\mathcal{X}$ and $\Pi_\mathcal{Y}$ the projections on $\mathcal{X}$ and $\mathcal{Y}$, respectively. Also, given two measurable spaces $\mathcal{E}$ and $\mathcal{F}$, a measurable map $f: \mathcal{E} \to \mathcal{F}$ and a measure $\mu$ on $\mathcal{E}$, we define the measure $f^*\mu$ on $\mathcal{F}$ in the natural way by $f^*\mu = \mu \circ f^{-1}$. We denote by $\delta_x$ the usual delta measure located at $x \in \mathcal{X}$. We also denote by $\mathcal{M}_1(\mathcal{X})$ and $\mathcal{M}_+(\mathcal{X})$ the set of probability measures and positive finite measures on $\mathcal{X}$, respectively. We endow both sets with the topology of weak convergence. If $\mathcal{X}$ is a Polish space, then we denote by $\mathcal{C}([0,T], \mathcal{X})$ the space of continuous functions $f: [0,T] \to \mathcal{X}$. We endow this space with the usual topology of uniform convergence.

We first define the structure of the class of noise processes that we are going to work with.

DEFINITION 2.1. A quadruple $(\mathcal{W}, \{\mathcal{P}_t\}_{t \geq 0}, \mathbf{P}_\omega, \{\theta_t\}_{t \geq 0})$ is called a *stationary noise process* if it satisfies the following:

(i) $\mathcal{W}$ is a Polish space;
(ii) $\mathcal{P}_t$ is a Feller transition semigroup on $\mathcal{W}$ which accepts $\mathbf{P}_\omega$ as its unique invariant measure;
(iii) the family $\{\theta_t\}_{t > 0}$ is a semiflow of measurable maps on $\mathcal{W}$ satisfying the property $\theta_t^* \mathcal{P}_t(x, \cdot) = \delta_x$ for every $x \in \mathcal{W}$ and every $t > 0$.

The following definition is taken from [7] and provides the general framework in which we are going to address the question of ergodicity.

DEFINITION 2.2. A *continuous stochastic dynamical system* (SDS) on the Polish space $\mathcal{X}$ over the stationary noise process $(\mathcal{W}, \{\mathcal{P}_t\}_{t \geq 0}, \mathbf{P}_\omega, \{\theta_t\}_{t \geq 0})$ is a map

$$\Lambda: \mathbf{R}_+ \times \mathcal{X} \times \mathcal{W} \to \mathcal{X}, \qquad (t, x, w) \mapsto \Lambda_t(x, w),$$

with the following properties.

(i) *Regularity of paths*: For every $T > 0$, $x \in \mathcal{X}$ and $w \in \mathcal{W}$, the map $\Phi_T(x, w) : [0, T] \to \mathcal{X}$ defined by

$$\Phi_T(x, w)(t) = \Lambda_t(x, \theta_{T-t} w)$$

belongs to $\mathcal{C}([0, T], \mathcal{X})$.

(ii) *Continuous dependence*: The map $(x, w) \mapsto \Phi_T(x, w)$ is continuous from $\mathcal{X} \times \mathcal{W}$ to $\mathcal{C}([0, T], \mathcal{X})$ for every $T > 0$.

(iii) *Cocycle property*: The family of maps $\Lambda_t$ satisfies

$$\begin{aligned}\Lambda_0(x, w) &= x, \\ \Lambda_{s+t}(x, w) &= \Lambda_s(\Lambda_t(x, \theta_s w), w),\end{aligned} \tag{2.1}$$

for all $s, t > 0$, all $x \in \mathcal{X}$ and all $w \in \mathcal{W}$.

Given an SDS as in Definition 2.2 and an initial condition $x_0 \in \mathcal{X}$, we now show how to use it to construct in a natural way an $\mathcal{X}$-valued stochastic process with initial condition $x_0$. First, given $t \geq 0$ and $(x, w) \in \mathcal{X} \times \mathcal{W}$, we construct a probability measure $\mathcal{Q}_t(x, w; \cdot)$ on $\mathcal{X} \times \mathcal{W}$ by

$$\mathcal{Q}_t(x, w; A \times B) = \int_B \delta_{\Lambda_t(x, w')}(A) \mathcal{P}_t(w, dw'), \tag{2.2}$$

where $\delta_x$ denotes the delta measure located at $x$, $A$ is a measurable subset of $\mathcal{X}$ and $B$ is a measurable subset of $\mathcal{W}$. One can show [7], Lemma 2.12, that the family of measures $\mathcal{Q}_t(x, w; \cdot)$ actually forms a Feller transition semigroup on $\mathcal{X} \times \mathcal{W}$ and if a probability measure $\mu$ on $\mathcal{X} \times \mathcal{W}$ satisfies $\Pi^*_\mathcal{W} \mu = \mathbf{P}_\omega$, then $\Pi^*_\mathcal{W} \mathcal{Q}_t \mu = \mathbf{P}_\omega$ for all $t > 0$. This suggests the following definition.

DEFINITION 2.3. Let $\Lambda$ be an SDS as above. Then a probability measure $\mu$ on $\mathcal{X} \times \mathcal{W}$ is called a *generalized initial condition* for $\Lambda$ if $\Pi^*_\mathcal{W} \mu = \mathbf{P}_\omega$. We denote by $\mathcal{M}_\Lambda$ the space of generalized initial conditions endowed with the topology of weak convergence. Elements of $\mathcal{M}_\Lambda$ of the form $\mu = \delta_x \times \mathbf{P}_\omega$ for some $x \in \mathcal{X}$ will be called *initial conditions*.

Given a generalized initial condition $\mu$, we construct a stochastic process $(x_t, w_t)$ on $\mathcal{X} \times \mathcal{W}$ by drawing its initial condition according to $\mu$ and then evolving it according to the transition semigroup $\mathcal{Q}_t$. The marginal $x_t$ of this process on $\mathcal{X}$ will be called the *process generated by $\Lambda$ for $\mu$*. We will denote by $\bar{\mathcal{Q}} \mu$ the law of this process [i.e., $\bar{\mathcal{Q}} \mu$ is a measure on $\mathcal{C}(\mathbf{R}_+, \mathcal{X})$].

DEFINITION 2.4. Two generalized initial conditions $\mu$ and $\nu$ of an SDS $\Lambda$ are *equivalent* if the processes generated by $\mu$ and $\nu$ are equal in law. In short, $\mu \sim \nu \Leftrightarrow \bar{\mathcal{Q}} \mu = \bar{\mathcal{Q}} \nu$.

We say that a generalized initial condition $\mu$ is *invariant* for the SDS $\Lambda$ if it is invariant for the Markov transition semigroup $\mathcal{Q}_t$ generated by $\Lambda$. Similarly, we call it *ergodic* if it is ergodic for $\mathcal{Q}_t$, that is, if the law of the stationary Markov process on $\mathcal{X} \times \mathcal{W}$ with transition probabilities $\mathcal{Q}_t$ and fixed-time marginal $\mu$ is ergodic for the shift map.

The following remark turns out to be very useful for the approach taken in this work.

LEMMA 2.5. *The map $\bar{\mathcal{Q}}$ preserves ergodicity in the sense that if $\mu \in \mathcal{M}_\Lambda$ is an ergodic invariant measure for the SDS $\Lambda$, then $\bar{\mathcal{Q}}\mu$ is an ergodic invariant measure for the shift map on $\mathcal{C}(\mathbf{R}_+, \mathcal{X})$.*

PROOF. This is an immediate consequence of the general fact that if $T$ and $\tilde{T}$ are two measurable transformations on measure spaces $\mathcal{E}$ and $\tilde{\mathcal{E}}$ and there exists a measurable map $f : \mathcal{E} \to \tilde{\mathcal{E}}$ such that $f \circ T = \tilde{T} \circ f$, then if a measure $\mu$ is ergodic for $T$, $f^*\mu$ is ergodic for $\tilde{T}$. □

REMARK 2.6. As a consequence of the above result, if $\mu, \nu \in \mathcal{M}_\Lambda$ are two ergodic invariant measures for the semigroup $\mathcal{Q}_t$, then either $\bar{\mathcal{Q}}\mu = \bar{\mathcal{Q}}\nu$ or $\bar{\mathcal{Q}}\mu$ and $\bar{\mathcal{Q}}\nu$ are mutually singular.

**3. An abstract ergodicity result.** The main motivation of this section is provided by the following well-known facts from the theory of Markov processes. Recall that a Markov process on a Polish space $\mathcal{X}$ with transition probabilities $\bar{\mathcal{P}}_t$ is called *topologically irreducible at time $t$* if $\bar{\mathcal{P}}_t(x, A) > 0$ for every $x \in \mathcal{X}$ and every open set $A \subset \mathcal{X}$. We call it simply *topologically irreducible* if there exists such a time.

It is called *strong Feller at time $t$* if $\bar{\mathcal{P}}_t\psi$ is continuous for every bounded measurable function $\psi : \mathcal{X} \to \mathbf{R}$. Here, we abused notation and again used the symbol $\bar{\mathcal{P}}_t$ to denote the corresponding semigroup acting on observables. It is immediate that the strong Feller property is equivalent to the continuity of the function $x \mapsto \bar{\mathcal{P}}_t(x, \cdot)$ if the space of probability measures on $\mathcal{X}$ is equipped with the topology of strong convergence. A standard result often attributed to Doob and Khas'minskii states the following.

THEOREM 3.1 (Doob–Khas'minskii). *If a Markov process on a Polish space $\mathcal{X}$ with transition probabilities $\bar{\mathcal{P}}_t$ is topologically irreducible and strong Feller, then it can have at most one invariant probability measure.*

In this section we introduce the *strong Feller* property and *irreducibility* in the abstract framework of SDS as laid out in the previous section. As already pointed out in [9], the strong Feller property as stated above is actually not easily amenable to generalization, mainly because the topology of strong convergence of measures is not metrizable. Instead of generalizing the notion of continuity of the transition probabilities in the topology of strong convergence, we will thus follow the approach laid out in [9] and provide a generalization of the notion of continuity in the *total variation* topology.

In this section, we consider, as before, a general SDS $\Lambda$ on a Polish space $\mathcal{X}$ with stationary noise process $(\mathcal{W}, \{\mathcal{P}_t\}_{t\geq 0}, \mathbf{P}_\omega, \{\theta_t\}_{t\geq 0})$. Remember that we introduced a linear map $\bar{\mathcal{Q}}$ from $\mathcal{M}_1(\mathcal{X} \times \mathcal{W})$ into $\mathcal{M}_1(\mathcal{C}(\mathbf{R}_+, \mathcal{X}))$ constructed as the law of the process on $\mathcal{X}$ with a given initial condition. Denoting by $R_t: \mathcal{C}(\mathbf{R}_+, \mathcal{X}) \to \mathcal{C}([t, \infty), \mathcal{X})$ the natural restriction map, we define the sets

$$\mathcal{N}_\mathcal{W}^t = \{(w, \tilde{w}) \in \mathcal{W}^2 | R_t^* \bar{\mathcal{Q}} \delta_{(x,w)} \sim R_t^* \bar{\mathcal{Q}} \delta_{(x,\tilde{w})} \, \forall x \in \mathcal{X}\},$$
$$\mathcal{N}_\mathcal{X}^t = \{(x, y, w) \in \mathcal{X}^2 \times \mathcal{W} | R_t^* \bar{\mathcal{Q}} \delta_{(x,w)} \not\perp R_t^* \bar{\mathcal{Q}} \delta_{(y,w)}\},$$
$$\mathcal{N}^t = \{(x, y, w, \tilde{w}) \in \mathcal{X}^2 \times \mathcal{W}^2 | (w, \tilde{w}) \in \mathcal{N}_\mathcal{W}^t \text{ and } (x, y, w) \in \mathcal{N}_\mathcal{X}^t\}.$$

Here and in the sequel, we write $\mu \sim \nu$ to denote that two measures $\mu$ and $\nu$ are mutually absolutely continuous and $\mu \perp \nu$ to denote that they are mutually singular. We will also use the notation $\mu \leq \nu$ as a shorthand for "$\mu(A) \leq \nu(A)$ for every measurable set $A$."

Note that, beside the symmetries obvious from the definitions, the set $\mathcal{N}^t$ has the property

$$(x, y, w, \tilde{w}) \in \mathcal{N}^t \to (x, y, \tilde{w}, w) \in \mathcal{N}^t.$$

Note, also, that in the Markovian case, $\bar{\mathcal{Q}} \delta_{(x,w)}$ is independent of $w$, so $\mathcal{N}_\mathcal{W}^t = \mathcal{W}^2$ and $\mathcal{N}^t$ can be considered as a subset of $\mathcal{X}^2$ for every $t > 0$.

Recall that a *coupling* between two measures $\mu$ and $\nu$ on a space $\mathcal{X}$ is a measure $\pi$ on $\mathcal{X}^2$ such that $\pi(A \times \mathcal{X}) = \mu(A)$ and $\pi(\mathcal{X} \times A) = \nu(A)$ for every measurable set $A \subset \mathcal{X}$. In the same spirit, we will say that $\pi$ is a *subcoupling* for $\mu$ and $\nu$ if $\pi(A \times \mathcal{X}) \leq \mu(A)$ and $\pi(\mathcal{X} \times A) \leq \nu(A)$.

Consider the map

$$\hat{\Lambda}_t: \mathcal{X}^2 \times \mathcal{W}^2 \to \mathcal{X}^2 \times \mathcal{W}^2,$$

defined as $\hat{\Lambda}_t(x, y, \omega, \tilde{\omega}) = (\Lambda_t(x, \omega), \Lambda_t(y, \tilde{\omega}), \omega, \tilde{\omega})$ for $(x, y) \in \mathcal{X}^2$ and $(\omega, \tilde{\omega}) \in \mathcal{W}^2$. We will abuse notation by also writing $\hat{\Lambda}(x, y)$ for the map from $\mathcal{W}^2$ to $\mathcal{X}^2 \times \mathcal{W}^2$ obtained by fixing the first two arguments.

With this notation in place, the abstract result laying the foundation for the present work is the following.

THEOREM 3.2. *Let $\hat{\Lambda}$ be as above and assume that there exists a time $t > 0$ and a jointly measurable map*

$$(x, y, w) \mapsto \mathcal{P}_t^{x,y}(w, \cdot) \in \mathcal{M}_+(\mathcal{W}^2)$$

*with the following properties:*

1. *the measure $\mathcal{P}_t^{x,y}(w,\cdot)$ is a subcoupling for $\mathcal{P}_t(w,\cdot)$ and $\mathcal{P}_t(w,\cdot)$ for every $(x,y,w) \in \mathcal{X}^2 \times \mathcal{W}$;*
2. *there exists $s > 0$ such that*

$$(3.1) \qquad (\hat{\Lambda}_t(x,y)^* \mathcal{P}_t^{x,y}(w,\cdot))(\mathcal{N}^s) > 0$$

*for every $(x,y,w) \in \mathcal{X}^2 \times \mathcal{W}$.*

*Then $\Lambda$ can have at most one invariant measure (up to the equivalence relation of Definition 2.4).*

PROOF. Assume, by contradiction, that $\mu$ and $\nu$ are two distinct ergodic invariant measures for the SDS $\Lambda$ such that $\bar{\mathcal{Q}}\mu \neq \bar{\mathcal{Q}}\nu$. We claim that there exist nonzero positive measures $\tilde{\mu}$, $\tilde{\nu}$ and $\bar{\nu}$ on $\mathcal{W} \times \mathcal{X}$ such that

$$R_s^* \bar{\mathcal{Q}} \mu \geq R_s^* \bar{\mathcal{Q}} \tilde{\mu} \not\perp R_s^* \bar{\mathcal{Q}} \bar{\nu} \sim R_s^* \bar{\mathcal{Q}} \tilde{\nu} \leq R_s^* \bar{\mathcal{Q}} \nu.$$

If we are able to construct such measures, it follows immediately that $R_s^* \bar{\mathcal{Q}}\mu$ and $R_s^* \bar{\mathcal{Q}}\nu$ are not mutually singular, thus leading to a contradiction, by Lemma 2.5. Let us consider the finite measure $\hat{\Lambda}_t(x,y)^* \mathcal{P}_t^{x,y}(w,\cdot)$ on $\mathcal{W}_1 \times \mathcal{W}_2 \times \mathcal{X}_1 \times \mathcal{X}_2$, where $(\mathcal{W}_1, \mathcal{W}_2)$ and $(\mathcal{X}_1, \mathcal{X}_2)$ denote two copies of $\mathcal{W}$ and $\mathcal{X}$, respectively. By assumption, there exist times $s > 0$ and $t > 0$ such that

$$(\hat{\Lambda}_t(x,y)^* \mathcal{P}_t^{x,y}(w,\cdot))(\mathcal{N}^s) > 0$$

for every $(x,y,w) \in \mathcal{X}^2 \times \mathcal{W}$. This shows that the measure $\theta^{(\mu,\nu)}$ on $\mathcal{N}^s$ defined by

$$\theta^{(\mu,\nu)}(A) := \int_{\mathcal{W}} \int_{\mathcal{X}^2} (\hat{\Lambda}_t(x,y)^* \mathcal{P}_t^{x,y}(w,\cdot))(A \cap \mathcal{N}^s) \mu(\omega, dx)\nu(\omega, dy) \mathbf{P}_w(d\omega),$$

is not identically 0. Here, $\mu(\omega, \cdot)$ and $\nu(\omega, \cdot)$ are the disintegrations of $\mu$ and $\nu$ with respect to $\mathbf{P}_w$. By using the hypothesis that $\mathcal{P}_t^{x,y}(w,\cdot)$ is a subcoupling for $\mathcal{P}_t(w,\cdot)$ and $\mathcal{P}_t(w,\cdot)$ for every $(x,y,w) \in \mathcal{X}^2 \times \mathcal{W}$, it follows immediately from the invariance of $\mu$ and $\nu$ that $\tilde{\mu} := \Pi^*_{\mathcal{W}_1 \times \mathcal{X}_1} \theta^{(\mu,\nu)}$ and $\tilde{\nu} := \Pi^*_{\mathcal{W}_2 \times \mathcal{X}_2} \theta^{(\mu,\nu)}$ are smaller than $\mu$ and $\nu$, respectively. Let us now consider the measure $\bar{\nu} := \Pi^*_{\mathcal{W}_1 \times \mathcal{X}_2} \theta^{(\mu,\nu)}$ on $\mathcal{W} \times \mathcal{X}$. The definitions of $\bar{\nu}$ and $\tilde{\nu}$ yield

$$R_s^* \bar{\mathcal{Q}} \bar{\nu} = \int_{\mathcal{N}^s} R_s^* \bar{\mathcal{Q}} \delta_{(y,\omega)} \theta^{(\mu,\nu)}(dx, dy, d\omega, d\tilde{\omega}),$$

$$R_s^* \bar{\mathcal{Q}} \tilde{\nu} = \int_{\mathcal{N}^s} R_s^* \bar{\mathcal{Q}} \delta_{(y,\tilde{\omega})} \theta^{(\mu,\nu)}(dx, dy, d\omega, d\tilde{\omega}).$$

Since $(\omega, \tilde{\omega}) \in \mathcal{N}^s_{\mathcal{W}}$, it follows that $R_s^* \bar{\mathcal{Q}} \bar{\nu} \sim R_s^* \bar{\mathcal{Q}} \tilde{\nu}$.

It remains to prove that $R_s^* \bar{\mathcal{Q}} \tilde{\mu} \not\perp R_s^* \bar{\mathcal{Q}} \bar{\nu}$. To see this, we observe that $\Upsilon := \Pi^*_{\mathcal{W}} \tilde{\mu} = \Pi^*_{\mathcal{W}} \bar{\nu}$ and therefore, by the triangle inequality and the fact that the measures $\tilde{\mu}$ and $\bar{\nu}$ give full measure to $\mathcal{N}^s$, one has the inequality

$$\|R_s^* \bar{\mathcal{Q}} \tilde{\mu} - R_s^* \bar{\mathcal{Q}} \bar{\nu}\|_{\mathrm{TV}}$$
$$\leq \int_{\mathcal{W}} \int_{\mathcal{X}^2} \|R_s^* \bar{\mathcal{Q}} \delta_{(x,\omega)} - R_s^* \bar{\mathcal{Q}} \delta_{(y,\omega)}\|_{\mathrm{TV}} \tilde{\mu}(\omega, dx) \bar{\nu}(\omega, dy) \Upsilon(d\omega)$$
$$< 2\Upsilon(\mathcal{W}) = \|R_s^* \bar{\mathcal{Q}} \tilde{\mu}\|_{\mathrm{TV}} + \|R_s^* \bar{\mathcal{Q}} \bar{\nu}\|_{\mathrm{TV}},$$

where $\tilde{\mu}(\omega, \cdot)$ and $\bar{\nu}(\omega, \cdot)$ are disintegrations of $\tilde{\mu}$ and $\bar{\nu}$, respectively, with respect to $\Upsilon$. The strict inequality from the first to the second line is an immediate consequence of the fact that the integral can be restricted to $\mathcal{N}_{\mathcal{X}}^s$ without changing its value. The claim $R_s^* \bar{\mathcal{Q}} \tilde{\mu} \not\perp R_s^* \bar{\mathcal{Q}} \bar{\nu}$ is then a consequence of the fact that if $\mu$ and $\nu$ are any two positive measures, then $\mu \perp \nu$ if and only if $\|\mu - \nu\|_{\mathrm{TV}} = \|\mu\|_{\mathrm{TV}} + \|\nu\|_{\mathrm{TV}}$.

Finally, note that since $\tilde{\mu} \leq \mu$ and $\tilde{\nu} \leq \nu$ by definition, we have $R_s^* \bar{\mathcal{Q}} \tilde{\mu} \leq R_s^* \bar{\mathcal{Q}} \mu$ and $R_s^* \bar{\mathcal{Q}} \tilde{\nu} \leq R_s^* \bar{\mathcal{Q}} \nu$. This completes the proof of the theorem. $\square$

The conditions of Theorem 3.2 do not appear to be easily verifiable at first sight. The remainder of this section is devoted to providing useful characterizations on the dynamics generated by an SDS $\Lambda$ on the state space $\mathcal{X}$ which give sufficient conditions for the assumptions in Theorem 3.2 to hold. It turns out that such properties are analogous to the strong Feller property and topological irreducibility in the Markovian setting.

DEFINITION 3.3. An SDS $\Lambda$ is said to be *strong Feller at time $t$* if there exists a jointly continuous function $\ell : \mathcal{X}^2 \times \mathcal{W} \to \mathbf{R}_+$ such that

(3.2) $$\|R_t^* \bar{\mathcal{Q}} \delta_{(x,\omega)} - R_t^* \bar{\mathcal{Q}} \delta_{(y,\omega)}\|_{\mathrm{TV}} \leq \ell(x, y, \omega)$$

and $\ell(x, x, \omega) = 0$ for every $x \in \mathcal{X}$ and every $\omega \in \mathcal{W}$.

REMARK 3.4. If the process is Markov in $\mathcal{X}$, then the total variation distance between $R_t^* \bar{\mathcal{Q}} \delta_{(x,\omega)}$ and $R_t^* \bar{\mathcal{Q}} \delta_{(y,\omega)}$ is equal to the total variation distance between the transition probabilities at time $t$ starting from $x$ and $y$, respectively. Therefore, Definition 3.3 reduces in this case to the statement "the transition probabilites at time $t$ are continuous in the total variation topology." This implies the usual strong Feller property but is not equivalent to it. However, it can be shown that if a Markov semigroup is strong Feller at time $t$, then the corresponding transition probabilities at time $2t$ are continuous in the total variation topology. This implies that for our purpose (where we are only interested behavior at large times anyway), Definition 3.3 is equivalent to the usual strong Feller property in the Markovian case.

DEFINITION 3.5. An SDS $\Lambda$ is said to be *topologically irreducible at time $t$* if for every $x \in \mathcal{X}$, $\omega \in \mathcal{W}$ and every nonempty open set $U \subset \mathcal{X}$, one has $\mathcal{Q}_t(x, \omega; \mathcal{W} \times U) > 0$.

REMARK 3.6. Since the dynamics which we are interested in take place in the state space $\mathcal{X}$, we do not generally require that the underlying Markov process generated by the semigroup $\mathcal{Q}_t$ be irreducible in the usual sense. In fact, the above definition is *much* weaker than irreducibility of the Markov semigroup $\mathcal{Q}_t$.

In the sequel, we will use the following notation. If $\mu$ is a finite measure on a measurable space $(Y, \mathcal{B}(Y))$ and $O \in \mathcal{B}(Y)$, then we write $\mu_{|O}(A) := \mu(A \cap O)$ for $A \in \mathcal{B}(Y)$. Next, we introduce a class of SDS which plays an important role in the theory.

DEFINITION 3.7. An SDS $\Lambda$ is said to be *quasi-Markovian* if for any two open sets $V, U$ in $\mathcal{W}$ and for every $t, s > 0$, there exists a measurable map $\omega \mapsto \mathcal{P}_s^{V,U}(\omega, \cdot) \in \mathcal{M}_+(\mathcal{W}^2)$ such that:

(i) the measure $\mathcal{P}_s^{V,U}(\omega, \cdot)$ is a subcoupling for $\mathcal{P}_s(\omega, \cdot)_{|V}$ and $\mathcal{P}_s(\omega, \cdot)_{|U}$ for every $\omega \in \mathcal{W}$,

(ii) one has $\mathcal{P}_s^{V,U}(\omega; \mathcal{N}_{\mathcal{W}}^t) > 0$ for every $\omega$ such that $\min\{\mathcal{P}_s(\omega; V), \mathcal{P}_s(\omega; U)\} > 0$.

REMARK 3.8. The terminology *quasi-Markovian* is motivated by the following fact. The process on $\mathcal{X}$ generated by the SDS $\Lambda$ is a Markov process in its own filtration precisely if $\bar{\mathcal{Q}}\delta_{(x,\omega)}$ is independent of $\omega$. In this case, $\mathcal{N}_{\mathcal{W}}^t = \mathcal{W}^2$ for every $t > 0$ and one can simply choose $\mathcal{P}_s^{V,U}(\omega; \cdot) = \mathcal{P}_s(\omega; \cdot)_{|V} \times \mathcal{P}_s(\omega; \cdot)_{|U}$ in the definition above so that $\Lambda$ is also quasi-Markovian.

REMARK 3.9. Definition 3.7 depends very weakly on the choice of the SDS $\Lambda$. It is weak in the sense that it does not take into account the fine details of the dynamics of $\Lambda$, but only how the noise enters the system. For example, we will show below that the solution to *any* SDE driven by fBm is always quasi-Markovian, without any restriction on the coefficients of the equation (1.1) other than what is required to obtain a well-posed SDS.

The last result of this section, which is the main abstract result of the present work, combines these definitions into a general criterion for an SDS to have a unique invariant measure. It is the analogue in our non-Markovian setting of Theorem 3.1 for Markovian systems.

THEOREM 3.10. *If there exist times $s > 0$ and $t > 0$ such that, a quasi-Markovian SDS $\Lambda$ is strong Feller at time $t$ and irreducible at time $s$, then satisfies the assumption of Theorem 3.2. In particular, it can have at most one invariant measure.*

PROOF. Since $\mathcal{W}$ is Polish, there exists a countable dense subset $\{w_n\}_{n\geq 0}$ and a metric $d_\mathcal{W}$ generating the topology of $\mathcal{W}$. Given this, we denote by $\{V_n\}_{n\geq 0}$ the countable collection of open balls with $1/2^k$ and center $w_m$ for all pairs of positive integers $k$ and $m$. We also choose a metric $d_\mathcal{X}$ on $\mathcal{X}$. We will henceforth denote by $\mathcal{B}_\rho^\mathcal{X}(x) \subset \mathcal{X}$ the open ball of radius $\rho$ centered at $x$ and similarly by $\mathcal{B}_\rho^\mathcal{W}(w) \subset \mathcal{W}$ the open balls in $\mathcal{W}$.

In order to verify the assumptions of Theorem 3.2, our aim is to find a measurable function

$$(3.3) \qquad (x, y, w) \mapsto (n_x, n_y)$$

and to define

$$(3.4) \qquad \mathcal{P}_s^{x,y}(\omega, \cdot) = \mathcal{P}_s^{V_{n_x}, V_{n_y}}(\omega, \cdot),$$

where the right-hand side uses the family of subcouplings from Definition 3.7.

Note, first, that it is possible to construct a measurable function $f : \mathcal{W} \to \mathcal{W}$ with the property that $f(\omega) \in \operatorname{supp} \mathcal{P}_s(\omega, \cdot)$ for every $\omega \in \mathcal{W}$. One way of constructing it is to define

$$n_1(\omega) = \inf\{n | \mathcal{P}_s(\omega, \mathcal{B}_1^\mathcal{W}(w_n)) > 0\}$$

and then, recursively,

$$n_k(\omega) = \inf\{n | \mathcal{P}_s(\omega, \mathcal{B}_{2^{-k}}^\mathcal{W}(w_n)) > 0 \text{ and } w_n \in \mathcal{B}_{2^{1-k}}(w_{n_{k-1}(\omega)})\}.$$

It follows from the density of the $w_k$ that $n_k(\omega) < \infty$ for every $k$ and every $\omega$. It follows from the construction that the sequence $w_{n_k(\omega)}$ is Cauchy for every $\omega \in \mathcal{W}$, so we can then set $f(w) = \lim_{k\to\infty} w_{n_k(\omega)}$. Since, by construction, $\mathcal{P}_s(\omega, \mathcal{B}_\rho^\mathcal{W}(f(\omega))) > 0$ for every $\rho > 0$, the function $f$ indeed has the required properties.

Define the map $\tilde{x} : \mathcal{X} \times \mathcal{W} \to \mathcal{W}$ by

$$\tilde{x}(x, \omega) = \Lambda_s(x, f(\omega)),$$

as well as the measurable map $r : \mathcal{X} \times \mathcal{W} \to \mathcal{W}$ by

$$(3.5) \quad r(x, \omega) = \sup\{\rho | \ell(x', \tilde{x}(x, \omega), \omega') \leq 1/4 \text{ for all } (x', \omega') \in \mathcal{B}_\rho(x, \omega)\},$$

where we use $\mathcal{B}_\rho(x, \omega)$ as shorthand for $\mathcal{B}_\rho^\mathcal{X}(\tilde{x}(x, \omega)) \times \mathcal{B}_\rho^\mathcal{W}(f(\omega))$. Since the function $\ell$ in Definition 3.3 is jointly continuous and vanishes on the diagonal, one has $r(x, \omega) > 0$ for every $x \in \mathcal{X}$ and every $\omega \in \mathcal{W}$.

Given these objects, we are now in a position to define $n_x$ and $n_y$. Consider $(x, y, \omega)$ to be given and set $\rho = r(x, \omega)$, $\tilde{\omega} = f(\omega)$ and $\tilde{x} = \Lambda_s(x, \tilde{\omega})$. We set

$$(3.6) \quad \begin{aligned} n_x &= \inf\{n | V_n \subset \mathcal{B}_\rho^\mathcal{W}(\tilde{\omega}), \ \tilde{\omega} \in V_n \text{ and } \Lambda_s(x, V_n) \subset \mathcal{B}_\rho^\mathcal{X}(\tilde{x})\}, \\ n_y &= \inf\{m | \mathcal{P}_s(\omega, V_m) > 0 \text{ and } \Lambda_s(y, V_m) \subset \mathcal{B}_\rho^\mathcal{X}(\tilde{x})\}. \end{aligned}$$

Since one has $\Lambda_s(x,\tilde{\omega}) = \tilde{x}$ by definition, it follows from the continuity of $\Lambda_s$ and the definition of the $V_n$ that $n_x$ is finite for every possible value of $x$ and $\omega$. The fact that $n_y$ is finite for every possible value of $x$, $y$ and $\omega$ follows from the topological irreducibility of $\Lambda$. This shows that $\mathcal{P}_s^{x,y}(\omega,\cdot)$ as given by (3.4) and (3.6) is a map satisfying the first assumption of Theorem 3.2.

It remains to show that (3.1) also holds. Since $\tilde{\omega} \in V_{n_x}$, one has $\mathcal{P}_s(\omega, V_{n_x}) > 0$. Furthermore $\mathcal{P}_s(\omega, V_{n_y}) > 0$, by construction, so $\mathcal{P}_s^{x,y}(\omega, \mathcal{N}_{\mathcal{W}}^t) > 0$. It follows from the definition of $\rho$ that one has

$$\|R_t^* \bar{\mathcal{Q}} \delta_{(\bar{x},w)} - R_t^* \bar{\mathcal{Q}} \delta_{(\bar{y},w)}\|_{\mathrm{TV}} \leq \tfrac{1}{2}$$

for every $(\bar{x}, \bar{y}) \in \mathcal{B}_\rho^{\mathcal{X}}(\tilde{x})^2$ and every $w \in V_{n_x}$. It follows immediately from the definition of $\mathcal{N}_t^{\mathcal{W}}$ that

(3.7) $$R_t^* \bar{\mathcal{Q}} \delta_{(\bar{x},w_x)} \not\perp R_t^* \bar{\mathcal{Q}} \delta_{(\bar{y},w_y)}$$

for every $(\bar{x}, \bar{y})$ as above, every $w_x \in V_{n_x}$, and every $w_y$ such that $(w_x, w_y) \in \mathcal{N}_t^{\mathcal{W}}$. Since one has $\Lambda(x, \omega_x) \in \mathcal{B}_\rho^{\mathcal{X}}(\tilde{x})$ and $\Lambda(y, \omega_y) \in \mathcal{B}_\rho^{\mathcal{X}}(\tilde{x})$ for every $(\omega_x, \omega_y) \in V_{n_x} \times V_{n_y}$, (3.1) is now a consequence of (3.7) and of condition (ii) in Definition 3.7. $\square$

The remainder of this article is devoted to showing that it is possible to associate an SDS to (1.1) in such a way that the assumptions of Theorem 3.10 are satisfied.

**4. Construction of the SDS.** In this section, we construct a continuous SDS induced by the SDE (1.1) in the sense that for every generalized initial condition $\mu$, the probability measure $\bar{\mathcal{Q}}\mu$ on the path space is an adapted solution to (1.1). This will then allow us to investigate ergodic properties of the SDE (1.1) according to the program laid out in the Introduction.

In this work, we will make use of the well-known Mandelbrot–Van Ness representation of the fBm [14]. The advantage of this representation is it is invariant under time-shifts, so it is natural for the study of ergodic properties. We may represent the two-sided fBm $B_H$ with Hurst parameter $H \in (0,1)$ in terms of a two-sided Brownian motion $B$ as

(4.1) $$B_H(t) = \alpha_H \int_{-\infty}^0 (-r)^{H-1/2}(dB(r+t) - dB(r))$$

for some $\alpha_H > 0$ (see [21] for more details).

4.1. *Noise space and the stationary noise process.* The aim of this section is to define the stationary noise process which will be used to investigate ergodic properties of the SDE (1.1) as laid out in Section 3. The main step is to consider a suitable Polish space in such a way that:

(1) there exists a Feller semigroup on the noise space which admits the fractional Brownian motion measure as its unique invariant measure;

(2) the solution map of the SDE (1.1) is continuous with respect to both the driving noise and initial conditions on $\mathbf{R}^d$.

One should note that such properties are closely related to the topology given on the noise space. In particular, if we consider white noise (or fractional noise with $H < 1/2$), property (2) could not be realized on any conventional space of paths, but one would have to work with rough paths instead [3].

The remainder of this section is devoted to choosing a topology which realizes (1) and (2). At first, one should note that by the properties of the fBm, it is convenient to work with Polish spaces defined by some norm which captures the Hölder continuity on bounded intervals and, at the same time, gives some kind of regularity at infinity. For this purpose, if $\gamma \in (0,1)$ and $\delta \in (0,1)$, then we define $\mathcal{W}_{(\gamma,\delta)}$ as the completion of $\mathcal{C}_0^\infty(\mathbf{R}_-;\mathbf{R})$ with respect to the norm

$$(4.2) \qquad \|\omega\|_{\gamma,\delta} = \sup_{t,s\in\mathbf{R}_-} \frac{|\omega(t) - \omega(s)|}{|t-s|^\gamma (1+|t|+|s|)^\delta}.$$

We write $\widetilde{\mathcal{W}}_{(\gamma,\delta)}$ for the corresponding space containing functions on the positive line instead of the negative one. We also write $\mathcal{W}_{(\gamma,\delta),T}$ and $\widetilde{\mathcal{W}}_{(\gamma,\delta),T}$ when we restrict the arguments to the intervals $[-T, 0]$ and $[0, T]$, respectively. It should be noted that $\|\cdot\|_{\gamma,\delta,T}$ is equivalent to the Hölderian norm $\|\cdot\|_\gamma$ for every $0 < T < \infty$. Moreover, $\mathcal{W}_{(\gamma,\delta)}$ is a separable Banach space. The following lemma states that there is a Wiener measure on $\mathcal{W}_{(\gamma,\delta)}$ for the fBm.

LEMMA 4.1. *Let $H \in (1/2, 1)$, $\gamma \in (1/2, H)$ and $\gamma + \delta \in (H, 1)$. There exists a Borel probability measure $\mathbf{P}_H$ on $\mathcal{W}_{(\gamma,\delta)}$ such that the canonical process associated to $\mathbf{P}_H$ is a fractional Brownian motion with Hurst parameter $H$.*

PROOF. One can show, as in [7], that the set of all continuous functions $w$ with $\|w\|_{\gamma',\delta'} < \infty$ for some $\gamma' > \gamma$ and some $\delta'$ such that $\delta' + \gamma' < \delta + \gamma$ is contained in $\mathcal{W}_{(\gamma,\delta)}$. The claim then follows from Kolmogorov's criterion and the behavior of fractional Brownian motion under time inversion. $\square$

One can similarly show that the two-sided fractional Brownian motion can be realized as the canonical process for some Borel measure $\tilde{\mathbf{P}}_H$ on $\mathcal{W}_{(\gamma,\delta)} \times \widetilde{\mathcal{W}}_{(\gamma,\delta)}$.

Consider now the operator $\mathcal{A}$ defined by

$$(4.3) \qquad \mathcal{A}\omega(t) = \beta_H \int_0^\infty \frac{1}{r} g\left(\frac{t}{r}\right) \omega(-r) \, dr,$$

where $g$ is given by

$$g(x) = x^{H-1/2} + (H - 3/2)x \int_0^1 \frac{(u+x)^{H-5/2}}{(1-u)^{H-1/2}} \, du$$

and $\beta_H = (H - 1/2)\alpha_H \alpha_{1-H}$. It is shown in Proposition A.2 below that $\mathcal{A}$ actually defines a bounded linear operator from $\mathcal{W}_{(\gamma,\delta)}$ into $\widetilde{\mathcal{W}}_{(\gamma,\delta)}$.

It will be convenient in the sequel to make use of fractional integration and differentiation. For $\alpha \in (0,1)$, we define the fractional integration operator $\mathcal{I}^\alpha$ and the corresponding fractional differentiation operator $\mathcal{D}^\alpha$ by

$$\mathcal{I}^\alpha f(t) = \frac{1}{\Gamma(\alpha)} \int_0^t (t-s)^{\alpha-1} f(s) \, ds,$$

(4.4)

$$\mathcal{D}^\alpha f(t) = \frac{1}{\Gamma(1-\alpha)} \frac{d}{dt} \int_0^t (t-s)^{-\alpha} f(s) \, ds.$$

For a comprehensive survey of the properties of these operators, see [20]. The most important property that we are going to use here is that $\mathcal{I}^\alpha$ and $\mathcal{D}^\alpha$ are each other's inverses. Furthermore, denoting by $\tau_h \colon w \mapsto w + h$ the shift map on $\widetilde{\mathcal{W}}_{(\gamma,\delta)}$, we have the following result.

LEMMA 4.2. *Let $\mathcal{H}(w, \cdot)$ be the transition kernel from $\mathcal{W}_{(\gamma,\delta)}$ to $\widetilde{\mathcal{W}}_{(\gamma,\delta)}$ given by*

$$\mathcal{H}(w, \cdot) = (\tau_{\mathcal{A}w} \circ \mathcal{I}^{H-1/2})^* \mathbf{W},$$

*where $\mathbf{W}$ is the Wiener measure over $\mathbf{R}_+$. Then $\mathcal{H}$ is the disintegration of $\tilde{\mathbf{P}}_H$ with respect to $\mathbf{P}_H$.*

PROOF. This is a lengthy, but straightforward, calculation, using the representation (4.1) for the fractional Brownian motion. □

We are now a in position to define our stationary noise process. For this, let us consider the one-sided Wiener shift $\theta_t \colon \mathcal{W}_{(\gamma,\delta)} \to \mathcal{W}_{(\gamma,\delta)}$ defined by

(4.5) $\qquad \theta_t \omega(s) := \omega(s - t) - \omega(-t), \qquad s \in \mathbf{R}_-, t \in \mathbf{R}_+.$

In order to construct the transition semigroup on $\mathcal{W}_{(\gamma,\delta)}$, we also introduce the "concatenation" function $M_t \colon \mathcal{W}_{(\gamma,\delta)} \times \widetilde{\mathcal{W}}_{(\gamma,\delta)} \to \mathcal{W}_{(\gamma,\delta)}$ defined by

(4.6) $\qquad M_t(\omega, \widetilde{\omega}) \stackrel{\text{def}}{=} \begin{cases} \widetilde{\omega}(s+t) - \widetilde{\omega}(t), & \text{if } -t < s, \\ \omega(s+t) - \widetilde{\omega}(t), & \text{if } s \leq -t \leq 0. \end{cases}$

With these definitions at hand, we set

(4.7) $\qquad \mathcal{P}_t(\omega, \cdot) := M_t^*(\omega, \cdot) \mathcal{H}(\omega, \cdot) \qquad \text{for } \omega \in \mathcal{W}_{(\gamma,\delta)} \text{ and } t \in \mathbf{R}_+.$

We are now in a position to state the following result.

LEMMA 4.3. *The quadruple $(\mathcal{W}_{(\gamma,\delta)}, \{\mathcal{P}_t\}, \mathbf{P}_H, \{\theta_t\})$ is a stationary noise process.*

PROOF. It is obvious from the definition that $M_t$ is continuous from $\mathcal{W}_{(\gamma,\delta)} \times \widetilde{\mathcal{W}}_{(\gamma,\delta)}$ to $\mathcal{W}_{(\gamma,\delta)}$. Moreover, the operator $\mathcal{A}$ is continuous from $\mathcal{W}_{(\gamma,\delta)}$ to $\widetilde{\mathcal{W}}_{(\gamma,\delta)}$. Therefore, we may conclude that $\mathcal{P}_t(\omega, \cdot)$ is a Feller semigroup on $\mathcal{W}_{(\gamma,\delta)}$. The fact that its only invariant measure is given by $\mathbf{P}_H$ is a straightforward calculation. All the other properties follow immediately from the definitions. $\square$

4.2. *Definition of the SDS and existence of an invariant measure.* Recall that we are concerned with the multidimensional SDE

$$(4.8) \quad X_t = X_0 + \int_0^t f(X_s)\,ds + \int_0^t \sigma(X_s)\,dB_H(t), \qquad 1/2 < H < 1,$$

where the integral with respect to $B_H$ is a pathwise Riemann–Stieltjes integral. This kind of equation has been studied by several authors (see, e.g., [3, 15, 18]) using different approaches, but we will mainly use the regularity results from [15].

Note that we actually need $d$ independent fractional Brownian motions to drive (4.8), so we consider $d$ copies of the stationary noise process defined in Lemma 4.3. With a slight abuse of notation, we again denote it by $(\mathcal{W}_{(\gamma,\delta)}, \{\mathcal{P}_t\}, \mathbf{P}_H, \{\theta_t\})$. We define the continuous shift operator $\mathcal{R}_T : \mathcal{W}_{(\gamma,\delta)} \to \widetilde{\mathcal{W}}_{(\gamma,\delta),T}$ by $(\mathcal{R}_T h)(t) := h(t-T) - h(-T)$ and set

$$(4.9) \qquad \Lambda : \mathbf{R}_+ \times \mathbf{R}^d \times \mathcal{W}_{(\gamma,\delta)} \to \mathbf{R}^d,$$

defined by $\Lambda_t(x, \omega) := \hat{\Phi}_t(x, \mathcal{R}_t \omega)(t)$, where $\hat{\Phi}_t : \mathbf{R}^d \times \widetilde{\mathcal{W}}_t \to \mathcal{C}([0,t], \mathbf{R}^d)$ is the solution map of equation (1.1) which depends on the initial conditions and the noise.

PROPOSITION 4.4. *Let $\Lambda$ be the function defined in (4.9). Then $\Lambda$ is a stochastic dynamical system over the stationary noise process $(\mathcal{W}_{(\gamma,\delta)}, \{\mathcal{P}_t\}, \mathbf{P}_H, \{\theta_t\})$. Moreover, for every generalized initial condition $\mu$, the process generated by $\Lambda$ for $\mu$ is an adapted weak solution of the SDE (1.1).*

PROOF. The regularity properties follow from [15], Theorem 5.1, and the fact that $\|\cdot\|_{\gamma,\delta,T}$ is equivalent to the Hölderian norm $\|\cdot\|_\gamma$ for every $0 < T < \infty$. The cocycle property is a direct consequence of the composition property of ODEs since we are dealing with pathwise solutions. Furthermore, it is obvious from Lemma 4.2 that every process generated by a generalized initial condition from $\Lambda$ is a weak solution of equation (4.8).

The adaptedness of the solution is a consequence of the construction since the transition probabilities $\mathcal{P}_t$ of the noise process do not depend on the solution $x$. □

To conclude this section, we study the problem of existence of the invariant measure for equation (4.8). The main difficulty in proving it comes from the pathwise stochastic integral. Suitable bounds on the stochastic integral, together with a dissipativity assumption, are sufficient to ensure existence. More specifically, in order to prove existence of the invariant measure, we make use of a Lyapunov function in the following sense.

DEFINITION 4.5. We say that $V: \mathbf{R}^d \to \mathbf{R}_+$ is a *Lyapunov function* for $\Lambda$ if $V^{-1}([0, K])$ is compact for every $0 \leq K < \infty$ and there exists a constant $C$ and a continuous function $\xi: [0, 1] \to \mathbf{R}_+$ with $\xi(1) < 1$ such that

$$\int V(x) \mathcal{Q}_t \mu(dx, dw) \leq C + \xi(t) \int V(x) \mu(dx, dw)$$

for every $t \in [0, 1]$ and every generalized initial condition $\mu$.

Note that this definition is slightly different from the one given in [7], but it is straightforward to check that the Krylov–Bogoliubov criterion nevertheless applies, so the existence of a Lyapunov function ensures the existence of an invariant measure.

PROPOSITION 4.6. *Assume that the hypotheses* (H1) *and* (H3) *hold. Then for every $p \geq 1$, the map $x \mapsto |x|^p$ is a Lyapunov function for the SDS $\Lambda$ defined above. Consequently, there exists at least one invariant measure for equation* (4.8) *and this invariant measure has moments of all orders.*

PROOF. The proof follows closely the proof of [15], Proposition 5.1, but we keep track on the dependence of the constants on the initial condition. Fix an arbitrary initial condition $x_0 \in \mathbf{R}^d$ and a realization $B_H$ of the fractional Brownian motion with Hurst parameter $H$. We define $x_t$ and $z_t$ on $t \in [0, 1]$ by

$$dx_t = f(x_t)\,dt + \sigma(x_t)\,dB_H(t), \qquad dz_t = f(z_t)\,dt,$$

where the initial condition for $x_t$ is are given by $x_0$ and the initial condition for $z_t$ is also given by $x_0$. We also define $y_t = x_t - z_t$ so that

$$y_t = \int_0^t (f(y_s + z_s) - f(z_s))\,ds + \int_0^t \sigma(y_s + z_s)\,dB_H(s) =: F_t + G_t.$$

Fix two arbitrary values $\alpha \in (1 - H, 1/2)$ and $\beta \in (1 - \alpha, H)$. Following [15], we define, for $t \in [0, 1]$,

$$h_t = |y_t| + \int_0^t \frac{|y_t - y_s|}{|t - s|^{\alpha+1}}\,ds$$

and verify that $h$ satisfies the conditions of the fractional Gronwall lemma [15], Lemma 7.6. Note, first, that the global Lipschitz continuity of $f$ implies that

$$|F_t| + \int_0^t \frac{|F_t - F_s|}{|t-s|^{\alpha+1}}\, ds \leq C \int_0^t |y_s|(t-s)^{-\alpha}\, ds.$$

Here and in the sequel, $C$ denotes a generic constant depending only on $\alpha$, $\beta$, $f$ and $\sigma$. In order to bound $G_t$, first note that since $\sigma$ is bounded,

$$|\sigma(z_s + y_s) - \sigma(z_r + y_r)| \leq C|z_s - z_r|^\beta + C|y_s - y_r|.$$

Also, note that the dissipativity condition on $f$ ensures that $|z_s - z_r| \leq C|s - r|(1 + |x_0|)$, so that

(4.10) $\quad |\sigma(z_s + y_s) - \sigma(z_r + y_r)| \leq C|s - r|^\beta(1 + |x_0|^\beta) + C|y_s - y_r|.$

It follows, in the same way as in [15], Proposition 5.1, that

$$|G_t| \leq C\|B_H\|_\beta \left( \int_0^t \frac{\sigma(z_s + y_s)}{s^\alpha}\, ds + \int_0^t \int_0^t \frac{|\sigma(z_s + y_s) - \sigma(z_r + y_r)|}{|s-r|^{\alpha+1}}\, dr\, ds \right)$$

$$\leq C\|B_H\|_\beta \left( 1 + |x_0|^\beta + \int_0^t \int_0^t \frac{|y_s - y_r|}{|s-r|^{\alpha+1}}\, dr\, ds \right).$$

One similarly obtains the bound

$$\int_0^t \frac{|G_t - G_s|}{|t-s|^{\alpha+1}}\, ds \leq C\|B_H\|_\beta \left( 1 + |x_0|^\beta + \int_0^t (t-s)^{-\alpha} \int_0^t \frac{|y_s - y_r|}{|s-r|^{\alpha+1}}\, dr\, ds \right).$$

Combining all of the above yields for $h$ that

$$|h_t| \leq C\|B_H\|_\beta \left( 1 + |x_0|^\beta + \int_0^t (1 + (t-s)^{-\alpha})h_s\, ds \right)$$

$$\leq C\|B_H\|_\beta \left( 1 + |x_0|^\beta + \int_0^t (t-s)^{-\alpha} t^\alpha s^{-\alpha} h_s\, ds \right).$$

The fractional Gronwall lemma [15], Lemma 7.6, then implies that

$$|h_t| \leq C\|B_H\|_\beta (1 + |x_0|^\beta) \exp(C\|B_H\|_\beta^{1/(1-\alpha)} t)$$

for every $t \in [0, 1]$. Furthermore, the dissipativity condition (H3) ensures the existence of $\gamma > 0$ such that $|z_t|^p \leq e^{-\gamma t}|x_0|^p + C$. Combining these bounds shows that for every $\eta > 0$, there exists a constant $C$ such that

$$|x_t|^p \leq (1+\eta)e^{-\gamma t}|x_0|^p + C\exp(C\|B_H\|_\beta^{1/(1-\alpha)} t).$$

Since $\|B_H\|_\beta$ is almost surely finite and $B_H$ is a Gaussian process, it has Gaussian tails by Fernique's theorem. This shows that there exists a constant $C$ such that, for every $t \in [0, 1]$, one has the bound

$$\int |x|^p (\mathcal{Q}_t \mu)(dx, dw) \leq (1+\eta)e^{-\gamma t} \int |x|^p \mu(dx, dw) + C,$$

uniformly over all generalized initial conditions $\mu$. Since $\eta$ was arbitrary and affects only the value of the constant $C$, one can choose it in such a way that $(1+\eta)e^{-\gamma} < 1$, thus concluding the proof of Proposition 4.6.  □

**5. Uniqueness of the invariant measure.** In order to simplify our notation, we fix once and for all $\gamma \in (1/2, H)$ and $\delta > 0$ such that $H < \gamma + \delta < 1$ and $\mathcal{W}$ we denote by the noise space with these indices. We also use the shorthand $\widetilde{\mathcal{W}}$ for $\widetilde{\mathcal{W}}_{(\gamma,\delta)}$ and $\widetilde{\mathcal{W}}_T$ for $\widetilde{\mathcal{W}}_{(\gamma,\delta),T}$.

The main goal of this section is to prove that the strong Feller property defined in Section 3 holds for the solutions to equation (4.8). This property, together with an irreducibility argument, will then provide the uniqueness of the invariant measure for our system. In the Markovian case, one efficient probabilistic tool to recover the strong Feller property is the Bismut–Elworthy–Li formula [5]. The main feature of this formula is that it provides bounds on the derivatives of a Markov semigroup which are independent of the bounds on the derivatives of the test function. In the non-Markovian case, one would expect to recover the strong Feller property by using a similar idea. In the language of the present article, given a measurable function $\varphi : \mathcal{C}([1,\infty); \mathbf{R}^d) \to \mathbf{R}$, a Bismut–Elworthy–Li type formula would be an expression for the derivative (in the $x$ variable) of the function $\bar{\mathcal{Q}}\varphi : \mathbf{R}^d \times \mathcal{W} \to \mathbf{R}$ defined by

$$(5.1) \qquad \bar{\mathcal{Q}}\varphi(x,w) := \int_{\mathcal{C}([1,\infty),\mathbf{R}^d)} \varphi(z) R_1^* \bar{\mathcal{Q}} \delta_{(x,w)}(dz)$$

that does not involve any derivative of $\varphi$.

The main technical difficulty one faces when trying to implement this program is that it seems to be very hard to prove that the Jacobian $J_{0,t}$ of the flow has bounded moments. We will overcome this by a cutoff procedure in Wiener space. The price we have to pay is that we are not able to show that $\bar{\mathcal{Q}}\varphi$ is differentiable in $x$, but only that it is continuous in $x$ for any given $w \in \mathcal{W}$ and that its modulus of continuity can be bounded by a function of $|x|$ and $\|w\|_{(\gamma,\delta)}$ only, uniformly over $\varphi$, with $\|\varphi\|_\infty \leq 1$. This is, however, sufficient for Definition 3.3 to apply.

We will need some basic lemmas concerning the smoothness of the solutions with respect to their initial conditions and with respect to the noise. We begin with an elementary regularity result. For sake of completeness, we give the details here. As in the previous section, we denote by

$$\hat{\Phi}_t : \mathbf{R}^d \times \widetilde{\mathcal{W}}_t \to \widetilde{\mathcal{W}}_t$$

the solution map for (1.1). (The fact that its image actually belongs to $\widetilde{\mathcal{W}}_t$ and not only to $\mathcal{C}([0,t], \mathbf{R}^d)$ is a consequence of Proposition 4.6 and of the regularity results from [15].) The main regularity result used in this section is the following.

LEMMA 5.1. *Assume that the coefficients of the SDE, $\sigma$ and $f$, satisfy hypotheses* (H1) *and* (H2). *Then the map $\hat{\Phi}_T$ is differentiable in both of its argument for each fixed $T > 0$.*

*Define the matrix-valued function $J_t = (D_x \Phi_T(x, \tilde{w}))(t)$. Then $J_t$ and its inverse $J_t^{-1}$, respectively, satisfy the equations*

$$(5.2) \qquad J_t = I + \int_0^t Df(x_s) J_s \, ds + \sum_{k=1}^d \int_0^t D\sigma_k(x_s) J_s \, d\tilde{w}_k(s),$$

$$(5.3) \qquad J_t^{-1} = I - \int_0^t J_s^{-1} Df(x_s) \, ds - \sum_{k=1}^d \int_0^t J_s^{-1} D\sigma_k(x_s) \, d\tilde{w}_k(s),$$

*where $I$ is the $d \times d$ identity matrix and where $x_s$ is shorthand for $(\Phi_T(x, \tilde{w}))(s)$. Here and below, we also defined $\sigma_k$ by $(\sigma_k)_i = \sigma_{k,i}$.*

*For any given $v \in \widetilde{\mathcal{W}}$, define the process $K_t^v = (D_{\tilde{w}} \Phi_T(x, \tilde{w}) v)(t)$. Then $K_t^v$ satisfies*

$$(5.4) \quad K_t^v = \int_0^t Df(x_s) K_s^v \, ds + \sum_{k=1}^d \int_0^t D\sigma_k(x_s) K_s^v \, d\tilde{w}_k(s) + \int_0^t \sigma(x_s) \, dv(s).$$

*In particular, defining $J_{s,t} = J_t J_s^{-1}$, the equation*

$$K_t^v = \int_0^t J_{s,t} \sigma(x_s) \, dv(s)$$

*holds. Furthermore, for every bounded set $A \subset \mathbf{R}^d$ of initial conditions and every $C > 0$, there exists $K > 0$ such that*

$$\|\Phi_T(x, \tilde{w})\|_\gamma + \|D_x \Phi_T(x, \tilde{w})\|_\gamma \leq K$$

*for every $x \in A$ and every $\tilde{w}$ with $\|\tilde{w}\|_\gamma \leq C$.*

PROOF. We know that $x_t = (\Phi_T(x, \tilde{w}))(t)$ is the unique solution to the equation

$$x_t = x + \int_0^t f(x_s) \, ds + \int_0^t \sigma(x_s) \, d\tilde{w}(s).$$

We can write this in the form $\Phi_T(x, \tilde{w}) = \mathcal{M}_T(x, \tilde{w}, \Phi_T(x, \tilde{w}))$, where $\mathcal{M}_T : \mathbf{R}^d \times \widetilde{\mathcal{W}}_T \times \widetilde{\mathcal{W}}_T \to \widetilde{\mathcal{W}}_T$ is a continuously differentiable map in all of its arguments.

Therefore, the claim follows from the implicit function theorem if we can show that the derivative of $\mathcal{M}_T$ in its last argument is of norm strictly smaller than 1. This is not true in general, but it holds for $T$ sufficiently small. The result then follows from the a priori bounds of Proposition 4.6, together with a standard gluing argument; see also [16] for a more detailed proof.

The equation for $J_t^{-1}$ is an immediate consequence of the chain rule. The last expression for $K_t^v$ is a consequence of the variation of constants formula. □

As usual with probabilistic proofs of regularization properties, the main ingredient in the proof of the strong Feller property will be an integration by parts formula. Since the point-of-view taken in most of this work is to study the solutions to (1.1) conditioned on a realization of the past of the driving fBm, the natural Gaussian space in which to perform this integration by parts will be the space $\widetilde{\mathcal{W}}$ endowed with the Gaussian measure $\mathcal{H}(0,\cdot)$. Note that it follows from (4.1) and the definition of $\mathcal{H}$ that the law of the canonical process $w$ under $\mathcal{H}(0,\cdot)$ is the same as the law of $\mathcal{I}^{H-1/2}w$ under the Wiener measure $\mathbf{W}$. Given any function $F$ on $\widetilde{\mathcal{W}}$, we can thus associate to it a function $\tilde{F} = F \circ \mathcal{I}^{H-1/2}$ on Wiener space. Note, also, that the reproducing kernel space $\mathcal{K}_H$ of $\mathcal{H}(0,\cdot)$ is precisely equal to those functions $v$ such that $\mathcal{D}^{H-1/2}v$ belongs to the reproducing kernel space $\mathcal{K}$ of $\mathbf{W}$. This allows us to carry over the whole formalism of Malliavin calculus [17].

Note that the properties of the fractional derivative and integral operators (4.4) are such that the notion of "adaptedness" with respect to the canonical process or the transformed process $\mathcal{I}^{H-1/2}$ agree, so a process $F_t$ is adapted to the increments of the canonical process if and only if $\tilde{F}_t$ is adapted to the increments of the canonical process. This allows us to speak of an adapted process in this setting without any ambiguity. In particular, the Malliavin derivative $\mathscr{D}F$ with respect to an adapted $\mathcal{K}_H$-valued process $v$ can be defined in a natural way by the equality

$$\langle \mathscr{D}F, v \rangle = \langle \mathscr{D}\tilde{F}, \mathcal{D}^{H-1/2}\tilde{v} \rangle \circ \mathcal{D}^{H-1/2}.$$

In particular, if $F$ is Fréchet differentiable with Fréchet derivative $DF$, it is also Malliavin differentiable and we have the equality

$$(5.5) \qquad \langle DF, v \rangle = \langle \mathscr{D}F, v \rangle,$$

where $\mathcal{K}_H$ is identified with a subspace of $\widetilde{\mathcal{W}}$ in the usual way.

The following bound is an immediate consequence of the integration by parts formula from Malliavin calculus, together with the representation (4.1) of fractional Brownian motion.

THEOREM 5.2. *Let $F, G: \widetilde{\mathcal{W}} \to \mathbf{R}$ be Malliavin differentiable functions such that $FG$, $F\|\mathscr{D}G\|_{\mathcal{K}_H}$ and $G\|\mathscr{D}F\|_{\mathcal{K}_H}$ are square integrable. Furthermore, let $v: \widetilde{\mathcal{W}} \to \widetilde{\mathcal{K}}_H$ be an adapted process such that $\|v\|_{\mathcal{K}_H}$ is square integrable. Then one has the bound*

$$(5.6) \quad |\mathbf{E}(\langle \mathscr{D}F, v \rangle G)| \leq (\mathbf{E}\|v\|_{\mathcal{K}_H}^2)^{1/2}((\mathbf{E}F^2 G^2)^{1/2} + (\mathbf{E}F^2 \|\mathscr{D}G\|_{\mathcal{K}_H}^2)^{1/2}).$$

*In the assumptions and the conclusion, expectations are taken with respect to* $\mathcal{H}(0, \cdot)$.

PROOF. Since the assumptions ensure that both sides of (5.6) are finite, a standard density argument shows that it suffices to check (5.6) under the assumption that $F$, $G$ and $v$ are in the space $\mathbf{D}^\infty$ of Malliavin smooth functions with bounded moments of all orders; see [17] for this notation.

Denoting by $\delta$ the adjoint of the derivative operator $\mathscr{D}$ in $\mathbf{L}^2(\widetilde{\mathcal{W}}, \mathcal{H}(0, \cdot))$, we then have

$$\mathbf{E}\langle \mathscr{D}FG, v \rangle = \mathbf{E}(FG\delta v).$$

On the other hand, one has $\mathscr{D}FG = G\mathscr{D}F + F\mathscr{D}G$, so

$$\mathbf{E}(\langle \mathscr{D}F, v\rangle G) = \mathbf{E}(FG\delta v) - \mathbf{E}(\langle \mathscr{D}G, v\rangle F).$$

It now suffices to note that the adaptedness of $v$ ensures that one can use the Itô isometry to get $\mathbf{E}|\delta v|^2 = \mathbf{E}\|v\|_{\mathcal{K}_H}^2$. $\square$

This result allows us to show the following.

PROPOSITION 5.3. *Assuming that* (H1)–(H3) *hold, the SDS constructed in the previous section has the strong Feller property of Definition* 3.3.

PROOF. Denote by $\widetilde{\mathcal{W}}_{[1,T]}$ the restriction of functions in $\widetilde{\mathcal{W}}$ to the interval $[1,T]$. For every bounded measurable function $\varphi_T : \widetilde{\mathcal{W}}_{[1,T]} \to \mathbf{R}$, define $\bar{\mathcal{Q}}\varphi : \mathbf{R}^d \times \mathcal{W} \to \mathbf{R}$ by (5.1). The strong Feller property (3.2) follows if we can show that

$$|\bar{\mathcal{Q}}\varphi_T(x, w) - \bar{\mathcal{Q}}\varphi_T(y, w)| \leq \ell(x, y, w),$$

uniformly over all $T > 1$ and all bounded Fréchet differentiable functions $\varphi$ with bounded Fréchet derivatives such that $\sup_X |\varphi(X)| \leq 1$.

Denoting by $\mathbf{E}_w$ for simplicity, expectations over $\widetilde{\mathcal{W}}$ with respect to the probability measure $\mathcal{H}(w, \cdot)$, we have

$$\bar{\mathcal{Q}}\varphi_T(x, w) = \mathbf{E}_w \varphi_T(\hat{\Phi}_T(x, \tilde{w})).$$

Setting $z_s = sx + (1-s)z$ and $\xi = x - y$, this yields

$$\begin{aligned}
&\bar{\mathcal{Q}}\varphi_T(y, w) - \bar{\mathcal{Q}}\varphi_T(x, w) \\
&\quad = \mathbf{E}_w \int_0^1 \langle D\varphi_T(\hat{\Phi}_T(z_s, \tilde{w})), D_x \hat{\Phi}_T(z_s, \tilde{w})\xi \rangle \, ds.
\end{aligned} \tag{5.7}$$

The problem at this point is that we do not have the a priori bounds on the Jacobian $D_x \hat{\Phi}_T$ that would be required in order to exchange the order of integration. This can, however, be overcome by the following cutoff procedure.

Note that we have the following result.

LEMMA 5.4. *The maps $N_T : \widetilde{\mathcal{W}} \to \mathbf{R}$ defined by*

$$N_T : \tilde{w} \mapsto \sup_{s \vee t \leq T} \frac{|\tilde{w}(t) - \tilde{w}(s)|}{|t-s|^\gamma}$$

*belong to $\mathbf{D}^{2,1}$.*

PROOF. Note that one actually has $N_T(\tilde{w}) = \sup_{s \vee t \leq T} \frac{\tilde{w}(t) - \tilde{w}(s)}{|t-s|^\gamma}$. The result then follows from [17], Proposition 2.1.3, and the fact ([20], Theorem 3.1) that $\mathcal{I}^{H-1/2}$ is bounded from $H^1$ to the space of $\gamma$-Hölder continuous functions. □

Now, let $\chi : \mathbf{R}_+ \to [0,1]$ be a smooth function such that $\chi(r) = 0$ for $r \geq 3$, $\chi(r) = 1$ for $r \leq 1$ and $|\chi'(r)| \leq 1$ for every $r$. Then the cutoff functions $\chi_{R,R'} : \tilde{w} \mapsto \chi(N_1(\tilde{w})/R)\chi(N_T(\tilde{w})/R')$ all belong to $\mathbf{D}^{2,1}$ and one has the following, obvious, bound result.

LEMMA 5.5. *There exists a constant $C > 0$ such that for every $T > 1$, there exists $R'_T$ such that the $\mathbf{D}^{2,1}$ norm of $\chi_{R,R'}$ is bounded by $C$ for every $R$ and every $R' \geq R'_T$.*

Denoting by $\delta \bar{\mathcal{Q}} \varphi^w_{x,y}$ (as short hand) the left-hand side of (5.7), we get the bound

$$\begin{aligned}|\delta \bar{\mathcal{Q}} \varphi^w_{x,y}| &\leq |\mathbf{E}_w((1 - \chi_{R,R'})(\varphi_T \circ \hat{\Phi}_T)(x, \cdot))| \\ &\quad + |\mathbf{E}_w((1 - \chi_{R,R'})(\varphi_T \circ \hat{\Phi}_T)(y, \cdot))| \\ &\quad + \left| \mathbf{E}_w \chi_{R,R'}(\tilde{w}) \int_0^1 \langle D\varphi_T(\hat{\Phi}_T(z_s, \tilde{w})), D_x \hat{\Phi}_T(z_s, \tilde{w})\xi \rangle \, ds \right| \\ &\stackrel{\text{def}}{=} T_1 + T_2 + T_3.\end{aligned}$$

Since $\varphi$ is bounded by 1, the first two terms in this expression are both bounded by

(5.8) $\quad T_1 + T_2 \leq 2\mathcal{H}(w, \{\tilde{w} | N_1(\tilde{w}) > R \text{ or } N_T(\tilde{w}) > R'\}).$

Concerning the last term, we can now use Lemma 5.1 to exchange the order of integration and obtain

$$T_3 = \left| \int_0^1 \mathbf{E}_w(\chi_{R,R'}(\tilde{w}) \langle D\varphi_T(\hat{\Phi}_T(z_s, \tilde{w})), D_x \hat{\Phi}_T(z_s, \tilde{w})\xi \rangle) \, ds \right|.$$

At this point, we use exactly the same trick as in the proof of the Bismut–Elworthy–Li formula [5] to transform the derivative with respect to $x$ into a Malliavin derivative with respect to the noise process $\tilde{w}$. Let $h:[0,1] \to \mathbf{R}$ be any smooth function with $\operatorname{supp} h \subset (0,1)$ and $\int_0^1 h(s)\,ds = 1$ and set

$$v_x(t) = \int_0^t h(s)\sigma^{-1}(x_s)J_s\xi\,ds,$$

where $x_s$ and $J_s$ are as in Lemma 5.1. It then follows from Lemma 5.1 that one has $K_t^v = J_t\xi$ for every $t \geq 1$. Therefore,

$$(D_x\Phi_T(x,\tilde{w})\xi)(t) = (D_w\Phi_T(x,\tilde{w})v_x)(t)$$

for every $t \geq 1$.

Since, on the other hand we assumed that $\varphi$ does not depend on the solution of the SDE up to time 1, this implies, by (5.5), that

$$(5.9) \qquad T_3 = \left| \int_0^1 \mathbf{E}_w(\chi_{R,R'}(\tilde{w})\langle \mathscr{D}(\varphi_T \circ \Phi_T)(z_s,\cdot), v_{z_s}\rangle)\,ds \right|.$$

Note, now, that it follows from Lemma 5.1 that $\frac{d}{dt}v_z(t)$ is $\gamma$-Hölder continuous with its norm bounded uniformly over $z$ in a ball of radius 1 around $x$ and over $N_1(\tilde{w}) \leq 3R$. Since $\chi_{R,R'}$ vanishes for $N_1(\tilde{w}) \geq 3R$, this shows that there exists a constant $C(R,x)$ depending on $R$ and $x$, but not on $R'$, such that we can replace $v$ in (5.9) by $\hat{v}$, defined by

$$\hat{v}(t) = \int_0^t \frac{dv}{ds}(s \wedge \tau)\,ds,$$

where $\tau$ is the stopping time defined as the first time that $\|\frac{d}{dt}v|_{[0,\tau]}\|_\gamma$ is greater or equal to $C(R,x)$. This ensures that $\|\frac{d}{dt}\hat{v}\|_\gamma \leq C(R,x)$ almost surely, while still being adapted.

In order to apply Theorem 5.2, it thus suffices to note that the fact that $\frac{d}{dt}\hat{v}_{z_s}(t) = 0$ for $t \geq 1$ implies that there exists a constant $C$ such that

$$\|\hat{v}\|_{\mathcal{K}_H}^2 = \int_0^\infty (\mathcal{D}^{H+1/2}\hat{v}_x(t))^2\,dt \leq C\left\|\frac{d}{dt}\hat{v}_x\right\|_\gamma.$$

Again using the fact that $\varphi$ is bounded by 1, we get, for $T_3$, the bound

$$T_3 \leq C\int_0^1 \mathbf{E}_w\left(\left\|\frac{d}{dt}\hat{v}_{z_s}\right\|_\gamma^2\right)^{1/2}(1 + (\mathbf{E}_w\|\mathscr{D}\chi_{R,R'}\|^2)^{1/2})\,ds$$
$$\leq C(R,x)\|x-y\|,$$

for all $y$ with $\|y-x\| \leq 1$. Here, we used Lemma 5.5 to obtain the last bound. Note that this bound does not depend on $R'$ and $T$, provided that $R'$ is larger than the value $R_0'$ found in Lemma 5.5. We can therefore let $R'$ go to $\infty$ and get

$$|\delta\bar{\mathcal{Q}}\varphi_{x,y}^w| \leq 2\mathcal{H}(w,\{\tilde{w}|N_1(\tilde{w}) > R\}) + C(R,x)\|x-y\|$$

for every $y$ with $\|y-x\| \leq 1$. Since both terms can be chosen to be continuous in $w$, $R$, $x$ and $y$ and since the first term tends to 0 as $R \to \infty$, the required bound follows at once. $\square$

REMARK 5.6. There is a direct relation between the integrability of $\|\sigma^{-1}(x_t)J_t\xi\|_\gamma$ and the continuity of $R_1^*\bar{\mathcal{Q}}\delta_{(x,\omega)}$ in the total variation topology. In fact, if $\|\sigma^{-1}(x_t)J_t\xi\|_\gamma$ has second moments, then $R_1^*\bar{\mathcal{Q}}\delta_{(x,\omega)}$ is not only continuous, but it is Lipschitz in $\mathbf{R}^d$. In fact, one then has the following generalization of the Bismut–Elworthy–Li formula:

$$(5.10) \quad D_\xi \bar{Q}\varphi(x,w) = \mathbf{E}_w\bigg((\varphi \circ \Phi)(x,\cdot)\int_0^\infty \mathcal{D}^{H-1/2}h(s)\sigma^{-1}(x_s)J_s\xi\, dB(s)\bigg),$$

where $B$ is the Brownian motion obtained from $\tilde{w}$ via $B = \mathcal{D}^{H-1/2}\tilde{w}$, $\bar{Q}\varphi$ is as in (5.1) and $h$ is a smooth function with support in $[0,1]$ and which integrates to 1. The main difficulty in getting good a priori bounds on the Jacobian comes from the fact that the diffusion coefficient for the SDE satisfied by the Jacobian is not globally bounded. One can check in [15] that, in general, the Jacobian has finite moments in $\gamma$-Hölder spaces for some $\gamma \in (1/2, H)$ if $H > 3/4$.

We conjecture that by using a Picard iteration, it should be possible to show integrability for the Jacobian in the supremum norm by realizing the pathwise Riemann–Stieltjes integral as a symmetric integral in the sense of Russo and Vallois [19]. In this case, there is a representation of the Riemann–Stieltjes integral in terms of the Skorokhod integral plus a trace term involving the Gross–Sobolev derivative; see [1] for more details. This may allow sufficient improvement of the existing estimates to get (5.10) directly, without requiring any cutoff procedure.

REMARK 5.7. Between the completion of this article and its publication, Hu and Nualart [10] announced bounds on the Jacobian for SDEs driven by fractional Brownian motion with $H > \frac{1}{2}$.

5.1. *Proof of the main result.* Similar to the Markovian case, the strong Feller property defined in Section 3 is not sufficient for uniqueness of invariant probability measures in the framework of SDS. In addition to the strong Feller property, we need an additional argument which provides the desired result of uniqueness. As discussed in Section 3, this will be achieved by means of an irreducibility argument jointly with the quasi-Markovian property for $\Lambda$. In general, irreducibility requires some kind of nondegeneracy of the diffusion term. As far as the quasi-Markovian property is concerned, it is a direct consequence of the properties of fractional Brownian motion. We first show that under (H1)–(H3), the SDS constructed above is topologically irreducible.

PROPOSITION 5.8. *Assume that the SDE* (4.8) *satisfies assumptions* (H1) *and* (H2). *Then the SDS* $\Lambda$ *induced by the SDE is irreducible at time* $t = 1$.

PROOF. Since everything is continuous, the proof of this is much easier than that of the original support theorem [22] and we can use the same argument as in [11], for example. The invertibility of $\sigma$ implies that the control system

$$\dot{x}_t = f(x_t) + \sigma(x_t)\dot{u}_t, \qquad t \in [0,1], x_0 \in \mathbf{R}^d,$$

is exactly controllable for every $x_0$. This shows that the solution map $\Phi_1 : \mathbf{R}^d \times \widetilde{\mathcal{W}} \to \mathbf{R}^d$ is surjective for every fixed value of the first argument. Since this map is also continuous, the claim follows from the fact that the topological support of $\mathcal{H}(w, \cdot)$ is all of $\widetilde{\mathcal{W}}$ (this is a consequence of the fact that it is a Gaussian measure whose Cameron–Martin space contains $\mathcal{C}_0^\infty$ and is thus a dense subspace of $\widetilde{\mathcal{W}}$). □

In the sequel, we will use the following notation. If $\mu_1$ and $\mu_2$ are positive measures with Radon–Nikodym derivatives $D_1$ and $D_2$, respectively, with respect to some common reference measure $\mu$, we define the measure $\mu_1 \wedge \mu_2$ by

$$(\mu_1 \wedge \mu_2)(dx) := \min\{D_1(x), D_2(x)\}\mu(dx).$$

Note that such a common measure $\mu$ can always be found (take $\mu = \mu_1 + \mu_2$, for example) and that the definition of $\mu_1 \wedge \mu_2$ does not depend on the choice of $\mu$.

Next, we prove that the SDS $\Lambda$ defined in (4.9) is quasi-Markovian over the stationary noise process $(\mathcal{W}, \{\mathcal{P}_t\}_{t\geq 0}, \mathbf{P}_H, \{\theta_t\}_{t\geq 0})$. The main technical estimate that we need for this is the following.

LEMMA 5.9. *Let $\mathcal{A}$ be the operator defined in Lemma* 4.2. *Then* $\mathcal{D}^{H+1/2}\mathcal{A}h \in \mathbf{L}^2(\mathbf{R}_+)$ *for every $h$ such that $h' \in \mathcal{C}_0^\infty(\mathbf{R}_-)$.*

PROOF. Note that a simple change of variables yields, for $\mathcal{A}h$, the formula

$$(5.11) \qquad (\mathcal{A}h)(t) = \int_0^\infty \frac{g(y)}{y}\bar{h}(t/y)\,dy,$$

where we set $\bar{h}(t) = h(-t)$ for convenience. This shows, in particular, that $\mathcal{A}h$ is smooth. Therefore, we have, for $\mathcal{D}^{H+1/2}\mathcal{A}h$, the expression

$$(5.12) \qquad \mathcal{D}^{H+1/2}\mathcal{A}h(t) = \int_0^\infty \frac{g(y)}{y^{H+3/2}}(\mathcal{D}^{H-1/2}h')(t/y)\,dy.$$

It follows immediately from the fact that $h' \in \mathcal{C}_0^\infty$ that there exists $C$ such that $\mathcal{D}^{H-1/2}h'$ is a smooth function bounded by $C\min\{1, t^{-H-1/2}\}$.

Therefore, one gets the bound

$$|\mathcal{D}^{H+1/2}\mathcal{A}h(t)| \leq C\left(t^{-H-1/2}\int_0^t \frac{g(y)}{y}\,dy + \int_t^\infty \frac{g(y)}{y^{H+3/2}}\,dy\right).$$

Using Lemma A.1, it is straightforward to get the bound

$$|\mathcal{D}^{H+1/2}\mathcal{A}h(t)| \leq C\min\{t^{-1}, t^{1/2-H}\},$$

where the constant $C$ depends on $h$. Since $t^{-1}$ is square integrable at $\infty$ and $t^{1/2-H}$ is square integrable at 0, this concludes the proof. $\square$

An immediate corollary of this is the following.

COROLLARY 5.10. *The set of pairs $(w, \tilde{w})$ such that $\tilde{w} - w \in \mathcal{C}_0^\infty(\mathbf{R}_-)$ is contained in $\mathcal{N}_\mathcal{W}^t$ for every $t > 0$.*

This eventually leads to the following result.

PROPOSITION 5.11. *The SDS in (4.9) induced by the SDE (4.8) is quasi-Markovian over the stationary noise process $(\mathcal{W}, \{\mathcal{P}_t\}_{t \geq 0}, \mathbf{P}_H, \{\theta_t\}_{t \geq 0})$ defined in Lemma 4.3.*

PROOF. Let us fix two nonempty open sets $U, V$ in $\mathcal{W}$ and two times, $s, t > 0$. As in the proof of Theorem 3.10, it is straightforward to construct measurable maps $f_U$ and $f_V$ from $\mathcal{W}$ to $\mathcal{W}$ with the property that $f_U(w) \in \operatorname{supp}\mathcal{P}_s(w, \cdot) \cap U$ for all $w$ such that $\mathcal{P}_s(w, U) > 0$, and similarly for $f_V$. Now, define a map $w \mapsto \varepsilon(w)$ by

$$\varepsilon(w) = \sup\{\varepsilon > 0 | B(f_U(w), \varepsilon) \subset U \text{ and } B(f_V(w), \varepsilon) \subset V\},$$

where we denote by $B(w, r)$ the ball of radius $r$ centered at $w$ in $\mathcal{W}$. Note that $\varepsilon(w) > 0$ on $A \stackrel{\text{def}}{=} \{w | \mathcal{P}_s(w, U) \wedge \mathcal{P}_s(w, V) > 0\}$.

Note, also, that the support of $\mathcal{P}_s(w, \cdot)$ consists precisely of those functions $\tilde{w}$ such that $\tilde{w}(t) - w(t+s)$ is constant for $t \leq -s$. Let $h_0(w) = f_V(w) - f_U(w)$ so that $h_0(w)$ is a function in $\mathcal{W}$ which is constant for times prior to $-s$. We now approximate $h_0$ by a smooth function $h$ with $h' \in \mathcal{C}_0^\infty$. This can, for example, be achieved by choosing two positive smooth functions $\psi$ and $\varphi$ such that $\psi(t) = 0$ for $t \notin [-2, -1]$ and $\int_{-1}^{-2}\psi(t)\,dt = 0$. Furthermore, $\varphi: \mathbf{R}_- \to [0, 1]$ is chosen to be decreasing and to satisfy $\varphi(-2) = 1$ and $\varphi(-1) = 0$. We then define

$$\mathcal{K}_\varepsilon h(t) = \varphi(t/\varepsilon)\int_{t-2\varepsilon}^{t-\varepsilon} h(r)\varepsilon^{-1}\psi((r-t)/\varepsilon)\,dr.$$

It is easy to check that $\mathcal{K}_\varepsilon h_0$ converges to $h_0$ strongly in $\mathcal{W}$ and that for every $\varepsilon > 0$, the derivative of $\mathcal{K}_\varepsilon h_0$ has support in $[-s+\varepsilon, -\varepsilon]$. We now define $h(w) = \mathcal{K}_{\delta(w)} h_0(w)$, where

$$\delta(w) = \sup\{\delta > 0 | \|\mathcal{K}_\delta h_0(w) - h_0(w)\| \leq \varepsilon(w)/2\}.$$

We thus constructed three measurable maps $w \mapsto \varepsilon(w)$, $w \mapsto f_U(w)$, and $w \mapsto h(w)$ such that, for every $w \in A$, the following properties hold:

(5.13a)
$$f_U(w) \in \operatorname{supp} \mathcal{P}_s(w, \cdot) \cap U,$$
$$B(f_U(w), \varepsilon(w)/2) \subset U,$$

(5.13b)
$$f_U(w) + h(w) \in \operatorname{supp} \mathcal{P}_s(w, \cdot) \cap V,$$
$$B(f_U(w) + h(w), \varepsilon(w)/2) \subset V.$$

Now, consider the maps $\mathcal{I}^1, \mathcal{I}^2 : \mathcal{W} \times \mathcal{W} \to \mathcal{W} \times \mathcal{W}$ defined by

$$\mathcal{I}^1(w, \tilde{w}) = (\tilde{w}, \tilde{w} + h(w)),$$
$$\mathcal{I}^2(w, \tilde{w}) = (\tilde{w} - h(w), \tilde{w}).$$

With this notation, we define

(5.14) $\quad \mathcal{P}_s^{U,V}(w, \cdot) = (\mathcal{I}^1(w, \cdot))^* \mathcal{P}_s(w, \cdot)_{|U} \wedge (\mathcal{I}^2(w, \cdot))^* \mathcal{P}_s(w, \cdot)_{|V}.$

It follows immediately from the definitions that $\mathcal{P}_s^{U,V}(w, \cdot)$ is a subcoupling for the measures $\mathcal{P}_s(w, \cdot)_{|U}$ and $\mathcal{P}_s(w, \cdot)_{|V}$. Since $h(w)$ has its derivative in $\mathcal{C}_0^\infty$, it is straightforward to check that it belongs to the reproducing kernel space of $\mathcal{P}_s(w, \cdot)$, so $(\mathcal{I}^1(w, \cdot))^* \mathcal{P}_s(w, \cdot)$ and $(\mathcal{I}^2(w, \cdot))^* \mathcal{P}_s(w, \cdot)$ are mutually absolutely continuous. This, together with (5.13), implies that $\mathcal{P}_s^{U,V}(w, \cdot) > 0$ for every $w \in A$, as required.

The fact that $\mathcal{P}_s^{U,V}(w, \mathcal{N}_\mathcal{W}^t) > 0$ is then an immediate consequence of Corollary 5.10. $\square$

Combining all of these results, we can now prove the main "concrete" result of this article.

THEOREM 5.12. *Under* (H1)–(H3), *there exists exactly one invariant probability measure for the SDS constructed in Section* 4.2.

PROOF. The existence of such an invariant measure is ensured by Proposition 4.6. Its uniqueness follows from Theorem 3.10, combined with Propositions 5.3, 5.8 and 5.11. $\square$

REMARK 5.13. Note that the solution to (1.1) obtained from the SDS constructed in Section 4.2 precisely coincide with the set of all adapted solutions to (1.1). Therefore, Theorem 1.1 is an immediate consequence of Theorem 5.12.

## APPENDIX

This section studies some of the properties of the operator $\mathcal{A}$ defined in (4.3). We first obtain the following bound.

LEMMA A.1. *Let $g: \mathbf{R}_+ \to \mathbf{R}$ be the function defined in Lemma 4.2. We then have $g(x) = \mathcal{O}(x)$ for $x \ll 1$ and $g(x) = \mathcal{O}(x^{H-1/2})$ for $x \gg 1$.*

PROOF. Since $g$ is smooth at every $x > 0$, we only need to check the result for $x \gg 1$ and $x \ll 1$. The behavior of $g(x)$ for $x \gg 1$ is straightforward since
$$x^{H-1/2} + (H - 3/2)x^{H-3/2} \le g(x) \le x^{H-1/2} \qquad \forall x \ge 0.$$

In order to treat the case $x \ll 1$, we rewrite $g$ as
$$g(x) = C_1 x (1+x)^{H-1/2} + C_2 x \int_0^1 (u+x)^{H-5/2}(1 - (1-u)^{1/2-H})\,du$$

for two constants $C_1$ and $C_2$. Note that for $x \ll 1$, the first term is $\mathcal{O}(x)$, so
$$|g(x)| \le Cx + Cx \int_0^1 u^{H-5/2}(1 - (1-u)^{1/2-H})\,du.$$

Now, note that $(1 - (1-u)^{1/2-H}) = \mathcal{O}(u)$ for $u \approx 0$ and that it diverges like $(1-u)^{1/2-H}$ for $u \approx 1$. Since, furthermore, $H > 1/2$, the function appearing under the integral is integrable, so $g(x) = \mathcal{O}(x)$. $\square$

This allows us to show the following.

PROPOSITION A.2. *Let $\gamma$ and $\delta$ be such that $1/2 < \gamma \le H$ and $H < \delta + \gamma < 1$. Then the operator $\mathcal{A}$ is bounded from $\mathcal{W}_{(\gamma,\delta)}$ into $\widetilde{\mathcal{W}}_{(\gamma,\delta)}$.*

PROOF. Fix $w \in \mathcal{W}_{(\gamma,\delta)}$ with $\|w\|_{(\gamma,\delta)} \le 1$ and consider two times $s$ and $t$ with $s < t$. Using (5.11), we obtain
$$|\mathcal{A}w(t) - \mathcal{A}w(s)| \le \int_0^\infty \frac{g(y)}{y} \frac{|t-s|^\gamma}{y^\gamma} \left(1 + \frac{t}{y}\right)^\delta dy,$$

so the claim follows if we can show that
$$\int_0^\infty \frac{g(y)}{y^{1+\gamma}} \left(1 + \frac{t}{y}\right)^\delta dy \le C(1+t)^\delta.$$

The left-hand side of this expression is bounded by
$$2 \int_0^t \frac{g(y)}{y^{1+\gamma}} \frac{t^\delta}{y^\delta}\,dy + 2 \int_t^\infty \frac{g(y)}{y^{1+\gamma}}\,dy \le 2t^\delta \int_0^\infty \frac{g(y)}{y^{1+\gamma+\delta}}\,dy + 2 \int_0^\infty \frac{g(y)}{y^{1+\gamma}}\,dy.$$

Now, note now that it follows immediately from Lemma A.1 that $g(y)y^{-\alpha}$ is integrable for every $\alpha \in (H + \frac{1}{2}, 2)$. This condition is satisfied for both $\alpha = 1 + \gamma$ and $\alpha = 1 + \gamma + \delta$, so the claim follows. $\square$

**Acknowledgments.** The second author would like to thank David Elworthy for his warm hospitality during a very pleasant and stimulating stay at the Mathematics Research Centre of the University of Warwick. The authors would also like to thank the referees for their careful reading of the manuscript.

MATHEMATICS INSTITUTE  
UNIVERSITY OF WARWICK  
COVENTRY, WARWICKSHIRE CV4 7AL  
UNITED KINGDOM  
E-MAIL: M.Hairer@warwick.ac.uk

DEPARTMENTO DE MATEMATICA  
UNICAMP  
6065, 13083-970 CAMPINAS  
SP, BRASIL  
E-MAIL: ohashi@ime.unicamp.br